\documentclass[12pt,english]{article}
\usepackage[T1]{fontenc}
\usepackage[latin9]{inputenc}
\usepackage{amsmath}
\usepackage{amssymb}

\makeatletter

\newtheorem{theorem}{Theorem}

\input{tcilatex}

\input{tcilatex}

\makeatother

\usepackage{babel}
\begin{document}
\begin{center}
{\Large{}The 3-D Spectrally-Hyperviscous Navier-Stokes Equations on
Bounded Domains with Zero Boundary Conditions}
\par\end{center}{\Large \par}

\begin{center}
{\Large{}\vspace{0.2in}
}
\par\end{center}{\Large \par}

\begin{center}
{\Large{}Joel Avrin }
\par\end{center}{\Large \par}

\begin{center}
{\Large{}Department of Mathematics and Statistics}
\par\end{center}{\Large \par}

\begin{center}
{\Large{}University of North Carolina at Charlotte}
\par\end{center}{\Large \par}

\begin{center}
{\Large{}jdavrin@uncc.edu \vspace{0.2in}
}
\par\end{center}{\Large \par}

\begin{center}
\underline{Abstract} 
\par\end{center}

We develop a mathematically and physically sound definition of the
spectrally-hyperviscous Navier-Stokes equations (SHNSE) on general
bounded domains $\Omega$ with zero (no-slip) boundary conditions
prescribed on $\varGamma=\partial\varOmega$. Previous successful
studies of the SHNSE have been limited to periodic-box domains, and
there are significant obtacles to overcome in extending the SHNSE
beyond this case; some of the numerical issues were discussed in {[}18{]},
and here we find a resolution of the theoretical issues. Beginning
with the basic hyperviscous case in which $A_{\varphi}u=(-\Delta)^{2}u$
is added to the NSE, we see that an additional boundary condition
is needed to make the operator $A_{\varphi}$ mathematically well-defined.
But this risks overdetermining the NSE system that the SHNSE is meant
to approximate, and indeed we show that the conditions $u=\frac{\partial u}{\partial\mathbf{n}}=0$
on $\varGamma$ as proposed in {[}28{]} are generally unphysical relative
to the NSE system as are the conditiions $u=\varDelta u=0$ on $\varGamma$. 

Taking an alternative approach we solve this issue successfully by
first applying the Leray projection $P$ to both sides of the NSE,
which requires making sense of the operator $P(-\Delta)^{2}$. Using
Helmholtz decomposition we show that $P(-\Delta)^{2}=A^{2}$ in $\Omega$
where $A=-P\Delta$ is the Stokes operator. The operator $A^{2}$
is well-defined and self-adjoint when equipped with the boundary conditions
$u=Au=0$ on $\varGamma$, and the fact that this formulation is physically
sound was recently shown in {[}6{]} in which it was shown that the
extra boundary condition $Au=0$ on $\varGamma$ necessarily holds
for the NSE if the forcing data is smooth enough. A version of the
SHNSE then results by setting $\mu PA_{\varphi}=\mu Q_{m}A^{2}u$
where $Q_{m}=I-P_{m}$ and $P_{m}$ is the projection onto the first
$m$ eigenspaces of $A$. Recent developments in {[}25{]}, {[}27{]}
give context to our results by clarifying the realizable impact of
SGS models.

With our new formulation of the SHNSE on general bounded domains in
hand we then establish foundational results, beginning with the existence
of globally regular solutions. Given that the SHNSE is meant to approximate
the NSE for small $\mu$ or large $m$, we establish this rigorously
in the general case by adapting the weak subsequence convergence results
of {[}5{]} to hold here. On intervals $[0,T]$ with a common $H^{1}$-bound
we deepen this sense of approximation by obtaining strong convergence.
First, by using estimates depending only on the common $H^{1}$-bound
to maximize computational applicability we show that SHNSE solutions
converge uniformly in $H^{1}$ to the NSE solution as either $\mu\rightarrow0$
or $m\rightarrow\infty$. Then in cases in which bootstrapped higher-order
bounds can be readily used we show that higher-order convergence results
hold. Our final results use the Stokes-pressure methodology developed
in {[}29{]}, {[}30{]} to recast the SHNSE in a form which like the
NSE reformulation in {[}29{]}, {[}30{]} is more adaptable to computation
and the specification of boundary values for the pressure.



\section{Introduction}

We consider the spectrally-hyperviscous Navier-Stokes equations (SHNSE)
for viscous incompressible homogeneous flow on a bounded domain $\Omega\subset\mathbb{R^{\textrm{3}}}$:
\begin{align}
u_{t}+\mu A_{\varphi}u-\nu\Delta u+\left(u\cdot\nabla\right)u+\nabla p=g,\tag{1.1a}\nonumber \\
\nabla\cdot u=0.\tag{1.1b}
\end{align}
Here $u=\left(u_{1},u_{2},u_{3}\right)$ is the velocity field of
the fluid, $g=\left(g_{1},g_{2},g_{3}\right)$ is the external force,
and $p$ is the pressure. We have that $u_{i}=u\left(x,t\right),g_{i}=g_{i}\left(x,t\right),i=1,2,3,$
and $p=p\left(x,t\right)$ where $x\in\Omega,$ a domain in $\mathbb{R}^{3}$.
In the treatments {[}3{]}, {[}4{]}, {[}5{]}, {[}18{]}, {[}19{]} of
(1.1) the domain $\Omega$ is assumed to be a periodic box, on which
after standardly ''moding out'' the constant vectors the operator
$B=-\Delta$ has eigenvalues $0<\lambda_{1}^{'}<\lambda_{2}^{'}<\cdots$
with corresponding eigenspaces $E_{1}^{'},E_{2}^{'},\cdots$. Let
$P_{m}^{'}$ be the projection onto $E_{1}^{'}\oplus\cdots\oplus E_{m}^{'}$,
let $Q_{m}^{'}=I-P_{m}^{'}$ and let $P_{E_{j}^{'}}^{'}$ be the projection
onto each $E_{j}^{'}$, then the basic general assumption on the operators
$A_{\varphi}$ considereded in {[}3{]}, {[}4{]}, {[}5{]} is that for
integers $\alpha\geq2$ we have that $A_{\varphi}\geq A_{m}\equiv Q_{m}^{'}A^{\alpha}$
in the sense of quadratic forms, i.e., $\left(A_{\varphi}v,v\right)\geq\left(A_{m}v,v\right)$
for smooth $v$. For purposes of regularity and for typical computational
purposes an applicable distinguished class (ADC) of the $A_{\varphi}$
identified in {[}3{]}, {[}4{]}, {[}5{]} and including the operators
considered in {[}18{]}, {[}19{]} satisfies $A_{\varphi}=\sum_{j=m_{0}+1}^{m}d_{j}(\lambda'_{j})^{\alpha}P'_{E'_{j}}+A_{m}$
where $0<m_{0}\leq m$ and $\{d_{j}\}_{j=m_{0}+1}^{m}$ is such that
$0<d_{j}\uparrow1$. Note that for $\mu=0$ the system (1.1) reduces
to the standard Navier-Stokes equations (NSE), and for $\mu>0$ we
have the hyperviscous Navier-Stokes equations (HNSE) in the special
case $m=0$ and $d_{j}\equiv1$. 

The motivation for the term $A_{\varphi}$ comes from the basic technique
of adding a subgrid-scale (SGS) stress tensor to the NSE to simulate
the dynamic effect of frequency scales too small to be resolved in
computations. Typically the SGS tensor is approximated by an extra
dissipative term, and the SGS configuration known as spectral-eddy
viscosity (SEV) was introduced in {[}23{]} to address some limitations
of these types of SGS models. In {[}10{]} it was suggested that SEV
could be approximated with a hyperviscous term, and the resulting
hyperviscous Navier-Stokes equations (HNSE) have been widely employed
computationally as in {[}7{]}, {[}8{]} and studied theoretically as
in {[}2{]}, {[}28{]}; see also the references contained therein. The
related spectral vanishing viscosity (SVV) method, introduced for
the study of gas dynamics in {[}35{]}, was applied to the incompressible
3-D NSE in {[}21{]}, the 2-D NSE in {[}34{]}, and to the 3-D NSE for
higher Reynolds numbers in e.g. {[}22{]}, {[}31{]}, {[}32{]}; see
in particular {[}31{]} in which Re = 768000. Hyperviscosity was implemented
spectrally in an application to conservation laws in {[}36{]}, and
the use of spectral hyperviscosity in application to the NSE was suggested
in {[}21{]}, discussed experimentally in {[}9{]}, and advocated in
{[}18{]}, {[}19{]}. The resulting spectrally-hyperviscous Navier-Stokes
equations (SHNSE) studied theoretically in {[}2{]}, {[}4{]}, {[}5{]},
{[}18{]} combine the SGS modeling and regularity of the HNSE with
the spectral accuracy philosophy of SVV; see e.g. {[}5{]}, {[}9{]},
{[}19{]} for further discussion of the SHNSE in relation to SGS, SEV,
and SVV modeling. Recent results in {[}25{]}, {[}27{]} showing that
fluctuations in turbulence do act dissipatively on the resolved scales
in the limit of time-averaging give further clarity regarding the
applicabilty of SGS modeling.

The SHNSE system has yet to reach its full potential as an SGS model
however, since the studies in {[}2{]}, {[}4{]}, {[}5{]}, {[}18{]}
were limited to periodic-box domains due to significant technical
issues. Some numerical issues were discussed in {[}19{]}, and to examine
the underlying theoretical issues we now consider (1.1) in the case
that $\Omega$ is a general bounded spatial domain in $\mathbb{R}^{3}$with
smooth (e.g. $C^{1}$) boundary $\Gamma$ on which we impose zero
(no-slip) boundary conditions. For simplicity of exposition we first
focus on the case $A_{\varphi}=B^{2}$. For $B^{2}$ to be well-defined
as a self-adjoint operator on $L^{2}(\Omega)$ an extra boundary condition
needs to be imposed. At the same time this extra condition needs to
be physically reasonable in the context of the NSE which the HNSE
and SHNSE are meant to approximate. In {[}28{]} the HNSE were configured
for bounded domains by assuming extra Neumann boundary conditions,
which in the case $A_{\varphi}=B^{2}$ are the conditions $u=\frac{\partial u}{\partial\mathbf{n}}=0$
on $\varGamma$. While mathematically well-defined, in practice when
types of Neumann conditions are imposed on the NSE the standard rigid-boundary
condition $u\cdot\mathbf{n}=0$ is retained but the tangental no-slip
condition generally is not (see e.g. {[}13{]}). Hence pairing the
no-slip condition with the extra Neumann condition is generally unphysical
(as borne out in conversation with a sample of applied practitioners)
or at best has severely limited application.

The boundary condition $u=\Delta u=0$ on $\Gamma$ seems mathematically
natural since $A_{\varphi}=(-\Delta)^{2}$ could be defined simply
as the square of the operator $B=-\Delta$ equipped with no-slip conditions.
But if $P$ denotes the Leray projection onto the solenoidal vectors,
then using the Helmholtz decomposition and the fact that $\nabla\cdot(-\Delta)u=(-\Delta)\nabla\cdot u=0$
in $\Omega$ for divergence-free $u$ we have as discussed in e.g.
{[}29{]}, {[}30{]} that the Stokes operator $A=-P\Delta$ satisfies
$Au=-\Delta u+\nabla p_{s}(u)$ where $p_{s}(u)$ solves the boundary-value
problem
\begin{align}
\Delta p_{s}(u)=0,x\in\Omega,\tag{1.2a}\nonumber \\
\mathbf{n}\cdot\nabla p_{s}(u)=\mathbf{n\cdot\Delta}u,x\in\Gamma.\tag{1.2b}
\end{align}
As noted in {[}29{]}, {[}30{]} the term $p_{s}(u)$ (referred to therein
as the Stokes pressure) satisfies $\nabla p_{s}(u)=(\Delta P-P\triangle)u=[\Delta,P]u$.
We see from (1.2) that if $\Delta u=0$ on $\Gamma$ then $p_{s}$
is a constant, hence $\nabla p_{s}=0$ and consequently $[\Delta,P]u=0$
which of course generally does not hold. Hence again we have an example
that has at best extremely limited physical applicability.

We now discuss how to successfully reconfigure (1.1) for general bounded
domains $\Omega$ in $\mathbb{R^{\textrm{3}}}$ with smooth boundary
$\partial\Omega=\Gamma$ and no-slip conditions $u=0$ imposed on
$\Gamma$. Our adaptation will preserve NSE physics and though we
assume a smoothness assumption on the forcing data we otherwise incur
no loss of generality. In particular we will use the results in {[}6{]}
which show that $A^{k}u=0$ on $\Gamma$ necessarily holds for solutions
$u$ of the NSE for any order $k$ provided that the forcing data
$f=Pg$ is in $D(A^{k-2})$. Hence for e.g. $\alpha=2$, while we
cannot assume that $u=\Delta u=0$ on the boundary, it necessarily
holds that $u=Au=0$ on $\Gamma$ for smooth enough $f$. But to use
these conditions we need a suitable reformulation of $A_{\varphi}$.
For the HNSE with e.g.. $\alpha=2$ we have that $A_{\varphi}=(-\Delta)^{2}$
in (1.1); applying $P$ to both sides as is standardly done for the
NSE requires making sense of the operator $P(-\Delta)^{2}$, and in
fact doing this is also useful toward our goal. With the decomposition
$P(-\Delta u)=-\Delta u+\nabla p_{s}(u)$ we obtain that 
\begin{align}
 & A^{2}u=P(-\Delta)P(-\Delta)\nonumber \\
 & =P(-\Delta)(-\Delta u+\nabla p_{s}(u))\nonumber \\
 & =P(-\Delta)^{2}u+P(-\Delta)(\nabla p_{s}(u))\tag{1.3}
\end{align}
and using (1.2a) and the commutivity of spatial derivatives inside
$\Omega$ we have that $P(-\Delta)(\nabla p_{s}(u))=-P\nabla(\Delta p_{s}(u))=0$
in $\Omega$. Combining with (1.3) we have in $\Omega$ that

\begin{equation}
A^{2}u=P(-\Delta)^{2}u.\tag{1.4}
\end{equation}
Since $A$ is of course well-defined assuming zero boundary conditions,
$P(-\Delta)^{2}u=A^{2}u$ is well-defined as a self-adjoint operator
assuming the conditions $u=Au=0$ on $\Gamma$. By induction using
(1.4) we have for any integer $\alpha\geq2$ that $P(-\Delta)^{\alpha}u=A^{\alpha}u$
is well-defined as a self-adjoint operator for smooth enough $f$
assuming the conditions $u=Au=\ldots=A^{\alpha-1}u=0$ on $\Gamma$. 

With (1.4) we thus we have a mathematically well-defined protocol
for defining the HNSE on general bounded domains $\Omega$ that is
by the results in {[}6{]} physically correct, assuming smooth enough
forcing data $f$ as above. Adding the term $\mu(-\Delta)^{\alpha}u$
to the NSE, applying $P$ to both sides, invoking (1.4), and associating
$A^{\alpha}$ with the boundary conditions $u=Au=\ldots=A^{\alpha-1}u=0$
on $\Gamma$ we obtain
\begin{align}
\frac{d}{dt}u+\mu A^{\alpha}u+\nu Au+P\left(u\cdot\nabla\right)u=f,\tag{1.5a}\nonumber \\
u(x,0)=u_{0}(x).\tag{1.5b}
\end{align}
From this new formulation of the HNSE we can derive a version of the
SHNSE for general bounded domains by using similar arguments as in
{[}5{]}. We let $0<\lambda_{1}<\lambda_{2}<\cdots$ represent the
eigenvalues of $A$ with corresponding eigenspaces $E_{1},E_{2},\cdots$.
Assuming as in {[}5{]} for large Reynolds numbers that $\mu$ is very
small, e.g. $\mu=\nu^{2}$ as in {[}7{]}, {[}8{]}, we envision similarly
as in {[}5{]} a cutoff $m$ such that $\mu\lambda_{j}^{\alpha}$ is
significant for $j\geq m$ and insignificant for $j\leq m$. Accordingly
we set $Q_{m}=I-P_{m}$ where $P_{m}$ is the projection onto $E_{1}\oplus\cdots\oplus E_{m}$,
and replace $A^{\alpha}$ in (1.5a) with operators $A_{\varphi}$
whose basic example is $A_{\varphi}=A_{m}\equiv Q_{m}A^{\alpha}$.
More generally, our basic assumtion on $A_{\varphi}$ is that $A_{\varphi}\geq A_{m}=Q_{m}A^{\alpha}$
in the sense of quadratic forms. This assumption is enough for many
of our results, and in particular is satisfied if $A_{\varphi}$ is
in the applicable distinguished class, with the ADC now defined using
the eigenvalues and eigenspaces of $A$. We assume in fact that $A_{\varphi}$
is in the ADC to simplify technical details of regularity in section
2; for example it follows straightforwardly if $A_{\varphi}$ is in
the ADC that the inequality $\left(A_{\varphi}u,A^{\theta}u\right)\geq\left(Q_{m}A^{\alpha+\theta}u,u\right)$
holds, which is a higher-order version of the quadratic-form inequality
$A_{\varphi}\geq A_{m}$.

With $A_{\varphi}$ redefined in this way we have the following formulation
of the SHNSE for general bounded domains
\begin{align}
\frac{d}{dt}u+\mu A_{\varphi}u+\nu Au+P\left(u\cdot\nabla\right)u=f,\tag{1.6a}\nonumber \\
u(x,0)=u_{0}(x).\tag{1.6b}
\end{align}
Applying $P_{n}$ to both sides of (1.6a) for $n\geq m$, for $f_{n}\equiv P_{n}f$
the Galerkin approximations to (1.6) are:
\begin{align}
\frac{d}{dt}u{}_{n}+\mu A_{\varphi}u_{n}+P_{n}P\left(u_{n}\cdot\nabla\right)u_{n}=f_{n},\tag{1.7a}\nonumber \\
u_{n}(x,0)=P_{n}u_{0}(x)\equiv u_{n,0}(x).\tag{1.7b}
\end{align}
Assuming for a constant $L$ that $\sup_{t\geq0}\left\Vert f\right\Vert _{2}\leq L$
and noting that $\left\Vert f_{n}\right\Vert _{2}\leq\left\Vert f\right\Vert _{2}$
and $\left\Vert u_{n,0}\right\Vert _{2}\leq\left\Vert u_{0}\right\Vert _{2}$
it follows straightforwardly for $v=u$ or $v=u_{n}$ that 
\begin{equation}
\left\Vert v\left(t\right)\right\Vert _{2}^{2}\leq\left\Vert u_{0}\right\Vert _{2}^{2}+\left(\frac{L}{\nu\lambda_{1}}\right)^{2}\equiv U_{L}^{2}.\tag{1.8}
\end{equation}

The development of (1.8) for solutions of (1.6) will be shown in section
2 below, and the arguments for solutions of (1.7) will follow similarly.
The estimate (1.8) is of course the same as the standard energy inequality
for the NSE; global existence and regularity for the ODE systems (1.7)
then follows by standard arguments. Bootstrapping from this we obtain
in section 2 a constant $U_{\mu,m}$ such that $\left\Vert A^{1/2}v(t)\right\Vert _{2}\leq U_{\mu,m}$
for all $t\geq0$, and from this regularity for (1.6) then follows
by slight modification of standard arguments for the NSE (see e.g.
{[}3{]}, {[}12{]}, {[}37{]}). Also in section 2 we will through modification
of arguments in {[}5{]} obtain weak subsequence convergence of solutions
of (1.6) to Leray solutions of the NSE as either $\mu\downarrow0$
or as $m_{0}\rightarrow\infty$, showing that in a suitable sense
the system (1.6) serves to approximate the NSE on general bounded
domains.

\begin{theorem} Let $\{u_{k}\}_{k=1}^{\infty}$ be the strong solutions
of (1.1) corresponding to $\mu=\mu_{k}$ where $\mu_{k}\downarrow0$
as $k\rightarrow\infty$. Then on each interval $[0,T]$ there exists
a subsequence, also denoted $\{u_{k}\}_{k=1}^{\infty}$, such that
$u_{k}\rightarrow v$ strongly in $L^{2}([0,T];H)$, $u_{k}\rightarrow v$
weakly in $L^{2}([0,T];V)$, $\frac{d}{dt}u_{k}\rightarrow\frac{d}{dt}v$
weakly in $L^{2}([0,T];PH^{-2})$, and $v$ is a Leray weak solution
of the NSE. If $\{u_{k}\}_{k=1}^{\infty}$ are the strong solutions
of (1.7) corresponding to $k=m_{0}$ and $k\rightarrow\infty$, then
on each interval $[0,T]$ there exists a subsequence, also denoted
$\{u_{k}\}_{k=1}^{\infty}$, and a $v\in L^{\infty}([0,T];H)\cap L^{2}([0,T];V)$
such that $u_{k}\rightarrow v$ strongly in $L^{2}([0,T];H)$, $u_{k}\rightarrow vu$
weakly in $L^{2}([0,T];V)$, $\frac{d}{dt}u_{k}\rightarrow\frac{d}{dt}v$
weakly in $L^{4/3}([0,T];V^{\prime})$, and $v$ is a Leray weak solution
of the NSE. \end{theorem} 

A more satisfying and robust sense of convergence would be established
by a full convergence result in a strong topology. This will be the
subject of our next results under the assumption of having an interval
of regularity $[0,T]$ for the NSE. Standard existence results for
the NSE (see e.g. the discussion and references in {[}12{]}, {[}15{]},
{[}26{]}, {[}37{]}) show that if $[0,T]$ is an interval over which
the $H^{1}$-norms of Galerkin solutions $v_{n}$ of the NSE are uniformly
bounded, then $v$ satisfies the same bound and can be continued as
a unique regular strong solution throughout $[0,T]$. On arbitrary
domains in 3-D such intervals of course include standard local intervals
of existence of strong solutions, as well as global intervals for
sufficiently small data. Results for larger data in 3-D have been
established on domains which are thin or have special symmetries (see
e.g. {[}1{]}, {[}24{]}, {[}26{]}, {[}33{]} and the references contained
therein). While it remains theoretically an open question whether
intervals of regularity can be constructed for arbitrarily large $T$
on arbitrary domains in 3-D with arbitrary data, the existence of
such intervals may in some cases be suggested by experimental observation
as noted e.g. in {[}20{]}. Such intervals $[0,T]$ are also intervals
over which a uniform $H^{1}$-bound exists independently of $\mu$
and $m$ for solutions of (1.1), which can be seen as follows: since
the semigroup $e^{-t\mu A_{\varphi}}$ is clearly a contraction on
the $H^{s}$-spaces and commutes with $A$, standard semigroup methods
show that any local interval of existence of strong solutions constructed
for the NSE with respect to the $H^{1}$-norm is also a local existence
interval for (1.1), and correspondingly any estimates used for the
NSE to continue such an interval with an $H^{1}$-bound can be used
similarly for (1.1). 

We thus assume with the above observations in mind that $[0,T]$ is
an interval over which the $H^{1}$-norm of the NSE solution $v$
and the $H^{1}$-norms of the solutions $u$ of (1.6) are uniformly
bounded by a common constant. Since by e.g. {[}17, Proposition 1.4{]}
the $H^{\gamma}$-norms and the $D(A^{\gamma/2})$-norms are equivalent
for all $v\in D(A^{\gamma/2})$, this common bound becomes a uniform
bound on the $D(A^{\gamma/2})$-norms, i.e. there exists a constant
$U_{1}$ such that for $h=u$ or $h=v$ 
\begin{equation}
\left\Vert A^{1/2}h(t)\right\Vert _{2}^{2}\leq U_{1}\tag{1.9}
\end{equation}
for all $t\in[0,T]$. 

On the interval $[0,T]$ for which (1.9) holds our next result obtains
strong $H^{1}$-convergence of solutions of (1.6) to strong NSE solutions
$v$. Despite the presence of higher-order terms represented by $A_{\varphi}$
the arguments depend only on (1.9) for greater utility of application
in computational settings. Iindeed it is more typical in such a setting
to assume and rely on no more than an $H^{1}$-bound, which in relation
to the NSE and SHNSE represents a finite energy in keeping with the
remarks above regarding intervals of regularity. The proof relies
significantly on new theoretical development that combines semigroup
methods with spectral decomposition techniques. 

\begin{theorem} Let $[0,T]$ be as above such that (1.9) holds, let
$u_{\mu}$ be the solutions of (1.6) for $\mu>0$, and let $u_{m}$
be the solutions of (1.6) for natural numbers $m$. Then for $f\in C[0,T];D(A^{1/2}))$
we have that $\sup_{0\leq t\leq T}\Vert A^{1/2}(u_{\mu}(t)-v(t))\Vert_{2}\rightarrow0$
as $\mu\downarrow0$ and that $\sup_{0\leq t\leq T}\Vert A^{1/2}(u_{m}(t)-v(t))\Vert_{2}\rightarrow0$
as $m\rightarrow\infty$ where $u$ is the unique global regular solution
of the NSE. \end{theorem}

Under additional smoothness assumptions on $f$ high-order bounds
can be bootstrapped from (1.7) (see Theorem 4 in section 2 below).
Though more unwieldy in computations, if we allow dependence on these
bounds our next result obtains high-order strong convergence. 

\begin{theorem} Let $[0,T]$ be as above such that (1.9) hoilds,
then for $f\in C[0,T];D(A^{(\theta+1)/2}))$ and solutions $u$ of
(1.6) and $v$ of the NSE we have for any natural numbers $\theta\geq2$
that $\sup_{0\leq t\leq T}\Vert A^{\theta/2}(u(t)-v(t))\Vert_{2}\rightarrow0$
as $\mu\downarrow0$ or as $m\rightarrow\infty$. \end{theorem}

We will prove Theorems 2 and 3 in section 3 below. In section 2 we
will make some preliminary observations, establish global regularity
as noted above, bootstrap higher-order bounds from (1.9), and sketch
the proof of Theorem 1. In section 4 we will again make use of the
Stokes pressure framework in {[}32{]}, {[}33{]} to derive reformulations
of (1.5) and (1.6) more amenable to computation and to the specification
of boundary values for the pressure. In section 5 we will make some
concluding remarks.

\section{Proof of Theorem 1}

We express the Sobolev inequalities on $\Omega$ in terms of the operator
$B=-\triangle:$
\begin{equation}
\left\Vert \upsilon\right\Vert _{q}\leq M_{1}\left\Vert B^{m_{1}/2}\upsilon\right\Vert _{2}^{1-\theta}\left\Vert B^{m_{2}/2}\upsilon\right\Vert _{2}^{\theta}\tag{2.1}
\end{equation}
for all $v\in D(B^{\theta/2})$ where $q\leq6/(3-2[(1-\theta)m_{1}+\theta m_{2}])$
and $M_{1}=M_{1}\left(\theta,q,m_{1},m_{2},\Omega\right)$. Here $B$
is equipped with zero boundary conditions as in the introduction.
By {[}17, Proposition 1.4{]} $D(A^{\gamma})$ is continuously embedded
into $H\cap H^{2\gamma}(\Omega)$ for any $\gamma\geq0$, and thus
we have for a constant $M_{0}=M_{0}\left(\theta,p,q,\Omega\right)$
and for $q$ as above that for all $v\in D(A^{\theta/2})$ 
\begin{equation}
\left\Vert \upsilon\right\Vert _{q}\leq M_{1}\left\Vert B^{m_{1}/2}\upsilon\right\Vert _{2}^{1-\theta}\left\Vert B^{m_{2}/2}\upsilon\right\Vert _{2}^{\theta}\leq M_{0}\left\Vert A^{m_{1}/2}\upsilon\right\Vert _{2}^{1-\theta}\left\Vert A^{m_{2}/2}\upsilon\right\Vert _{2}^{\theta}.\tag{2.2}
\end{equation}
Taking the inner product of both sides of (1.8a) with $Q_{k}Au_{n}$
for $m\leq k<n$ we have that 
\begin{equation}
\frac{1}{2}\frac{d}{dt}\left\Vert A^{1/2}Q_{k}u_{n}\right\Vert _{2}^{2}+\mu\Vert Q_{k}A^{3/2}u_{n}\Vert_{2}^{2}+\nu\left\Vert Q_{k}Au_{n}\right\Vert _{2}^{2}+(P_{n}P\left(u_{n}\cdot\nabla\right)u_{n},Q_{k}Au_{n})=(f_{n},Q_{k}Au_{n})\tag{2.3}
\end{equation}
 where we note that $(v,Q_{k}Au_{n})=(Q_{k}v,Q_{k}Au_{n})=(A^{1/2}Q_{k}v,Q_{k}A^{1/2}u_{n})=(A^{1/2}Q_{k}v,A^{1/2}Q_{k}u_{n})$
for any $v\in D(A^{1/2})$.

Now by the Cauchy-Schwartz and Young's inequalities 
\begin{equation}
|(f_{n},Q_{k}Au_{n})|\leq\frac{\nu}{2}\left\Vert Q_{k}Au_{n}\right\Vert _{2}^{2}+\frac{1}{2\nu}\left\Vert f\right\Vert _{2}^{2}\tag{2.4}
\end{equation}
 where we also use the fact that $\Vert P_{n}v\Vert_{2}\leq\Vert v\Vert_{2}$
for any $v\in H$. in similar fashion to the line after (2.3) we have
that $(P_{n}P\left(u_{n}\cdot\nabla\right)u_{n},Q_{k}Au_{n})=(P_{n}A^{-1/2}P\left(u_{n}\cdot\nabla\right)u_{n},Q_{k}A^{3/2}u_{n})$
so again by the Cauchy-Schwartz and Young's inequalities 
\begin{equation}
|(P_{n}P\left(u_{n}\cdot\nabla\right)u_{n},Q_{k}Au_{n})|\leq\frac{1}{2\mu}\left\Vert A^{-1/2}P\left(u_{n}\cdot\nabla\right)u_{n}\right\Vert _{2}^{2}+\frac{\mu}{2}\left\Vert ,Q_{k}A^{3/2}u_{n}\right\Vert _{2}^{2}.\tag{2.5}
\end{equation}
 Combining (2.4) and (2.5) with (2.3) and multiplying by 2 we have
that
\begin{align}
\frac{d}{dt}\left\Vert A^{1/2}Q_{k}u_{n}\right\Vert _{2}^{2}+\mu\Vert Q_{k}A^{3/2}u_{n}\Vert_{2}^{2}+\nu\left\Vert Q_{k}Au_{n}\right\Vert _{2}^{2}\nonumber \\
\leq\frac{1}{\mu}\left\Vert A^{-1/2}P\left(u_{n}\cdot\nabla\right)u_{n}\right\Vert _{2}^{2}+\frac{1}{\nu}\left\Vert f\right\Vert _{2}^{2}.\tag{2.6}
\end{align}
Now $A^{-1/2}P\left(u_{n}\cdot\nabla\right)u_{n}=A^{-1/2}Pdiv(u_{n}\otimes u_{n})$
for the appropriate tensor product $\otimes$ while $A^{-1/2}Pdiv\equiv T$
is a bounded operator on $H$ (see e.g. {[}17, Lemma 1.3{]}); thus
using (2.2) there is a constant $M_{1}$ such that
\begin{align}
\left\Vert A^{-1/2}P\left(u_{n}\cdot\nabla\right)u_{n}\right\Vert _{2}^{2}\leq\Vert T\Vert_{2}^{2}\left\Vert u_{n}\otimes u_{n}\right\Vert _{2}^{2}\nonumber \\
\leq\Vert T\Vert_{2}^{2}\parallel u_{n}\parallel_{2}^{2}\parallel u_{n}\parallel_{\infty}^{2}\leq M_{1}\Vert T\Vert_{2}^{2}\parallel u_{n}\parallel_{2}^{2}\parallel u_{n}\Vert_{2}\parallel A^{3/2}u_{n}\parallel_{2}.\tag{2.7}
\end{align}
Combining with (2.6) while neglecting the term $\nu\left\Vert Q_{k}Au_{n}\right\Vert _{2}^{2}$
on the left-hand side we have that
\begin{align}
\frac{d}{dt}\left\Vert A^{1/2}Q_{k}u_{n}\right\Vert _{2}^{2}+\mu\Vert Q_{k}A^{3/2}u_{n}\Vert_{2}^{2}\nonumber \\
\leq\frac{1}{\mu}\Vert T\Vert_{2}^{2}\parallel u_{n}\parallel_{2}^{2}M_{1}^{2}\parallel u_{n}\parallel_{2}\parallel A^{3/2}u_{n}\parallel_{2}+\frac{1}{\nu}\left\Vert f\right\Vert _{2}^{2}.\tag{2.8}
\end{align}
Combining terms in (2.8), setting $K_{1}=\Vert T\Vert_{2}^{2}M_{1}^{2}$
and using (1.11) we obtain that
\begin{align}
\frac{d}{dt}\left\Vert A^{1/2}Q_{k}u_{n}\right\Vert _{2}^{2}+\mu\Vert Q_{k}A^{3/2}u_{n}\Vert_{2}^{2}\nonumber \\
\leq(\frac{1}{\mu}K_{1}U_{L}^{3/2})\parallel A^{3/2}u_{n}\parallel_{2}+\frac{1}{\nu}L^{2}.\tag{2.9}
\end{align}
Applying Young's inequality to (2.9) we have that 
\begin{align}
 & \frac{d}{dt}\left\Vert A^{1/2}Q_{k}u_{n}\right\Vert _{2}^{2}+\mu\Vert Q_{k}A^{3/2}u_{n}\Vert_{2}^{2}\nonumber \\
 & \leq\frac{1}{2\mu^{2}}K_{1}U_{L}^{3}+\frac{\mu}{2}\parallel A^{3/2}u_{n}\parallel_{2}^{2}+\frac{1}{\nu}L^{2}\nonumber \\
 & =\frac{1}{2\mu^{2}}K_{1}U_{L}^{3}+\frac{\mu}{2}\parallel A^{3/2}P_{k}u_{n}\parallel_{2}^{2}+\frac{\mu}{2}\parallel A^{3/2}Q_{k}u_{n}\parallel_{2}^{2}+\frac{1}{\nu}L^{2}.\tag{2.10}
\end{align}
But $\parallel A^{3/2}P_{k}u_{n}\parallel_{2}^{2}\leq\lambda_{k}^{3}\parallel P_{k}u_{n}\parallel_{2}^{2}\leq\lambda_{k}^{3}\parallel u_{n}\parallel_{2}^{2}\leq\lambda_{k}^{3}U_{L}$
so using this in (2.10), subtracting, and using Poincaré's inequality
we obtain that
\begin{align}
\frac{d}{dt}\left\Vert A^{1/2}Q_{k}u_{n}\right\Vert _{2}^{2}+\frac{\mu}{2}\lambda_{k+1}^{2}\Vert Q_{k}A^{1/2}u_{n}\Vert_{2}^{2}\nonumber \\
\leq\frac{1}{2\mu^{2}}K_{1}U_{L}^{3}+\frac{\mu}{2}\lambda_{k}^{3}U_{L}+\frac{1}{\nu}L^{2}.\tag{2.11}
\end{align}
Integrating both sides of (2.11) we obtain for $d=\mu/2$ that 
\begin{equation}
\left\Vert A^{1/2}Q_{k}u_{n}\right\Vert _{2}^{2}\leq\left\Vert A^{1/2}Q_{k}u_{n,}{}_{0}\right\Vert _{2}^{2}e^{-d\lambda_{k+1}^{2}t}+\int_{0}^{T}(\frac{1}{2\mu^{2}}K_{1}U_{L}^{3}+\frac{\mu}{2}\lambda_{k}^{3}U_{L}+\frac{1}{\nu}L^{2})e^{-d\lambda_{k+1}^{2}(t-s)}ds\tag{2.12}
\end{equation}
from which we obtain in similar fashion to the development leading
to (1.11) that
\begin{equation}
\left\Vert A^{1/2}Q_{k}u_{n}(t)\right\Vert _{2}^{2}\leq\left\Vert A^{1/2}u_{0}\right\Vert _{2}^{2}+\frac{1}{d\lambda_{k+1}^{2}}(\frac{1}{2\mu^{2}}K_{1}U_{L}^{3}+\frac{\mu}{2}\lambda_{k}^{3}U_{L}+\frac{1}{\nu}L^{2}).\tag{2.13}
\end{equation}
But $\parallel A^{1/2}P_{k}u_{n}\parallel_{2}^{2}\leq\lambda_{k}\parallel P_{k}u_{n}\parallel_{2}^{2}\leq\lambda_{k}\parallel u_{n}\parallel_{2}^{2}\leq\lambda_{k}U_{L}$
so since $\parallel A^{1/2}u_{n}\parallel_{2}^{2}=\parallel A^{1/2}P_{k}u_{n}\parallel_{2}^{2}+\parallel A^{1/2}Q_{k}u_{n}\parallel_{2}^{2}$
we have, combining with (2.13), that
\begin{equation}
\left\Vert A^{1/2}u_{n}(t)\right\Vert _{2}^{2}\leq\lambda_{k}U_{L}+\left\Vert A^{1/2}u_{0}\right\Vert _{2}^{2}+\frac{1}{d\lambda_{k+1}^{2}}(\frac{1}{2\mu^{2}}K_{1}U_{L}^{3}+\frac{\mu}{2}\lambda_{k}^{3}U_{L}+\frac{1}{\nu}L^{2}).\tag{2.14}
\end{equation}
Through slight modification of these argumesnts we see that (2.14)
holds with $u_{n}$ replaced by solutuions $u$ of (1.6). Thus we
have for $v=u$ or $v=u_{n}$ that $\left\Vert A^{1/2}v(t)\right\Vert _{2}\leq U_{\mu,m}$
as in the remarks following (1.8) above with $U_{\mu,m}$ denoting
the right-hand side of (2.14).

We now show how the arguments in {[}5{]} can be modified to obtain
Theorem 1. Typically terms like $\left\Vert A^{-\beta}\left(u_{n}\cdot\nabla\right)u_{n}\right\Vert _{2}$
need to be estimated for various $\beta\geq0$, and here the corresponding
term is $\left\Vert A^{-\beta}P\left(u_{n}\cdot\nabla\right)u_{n}\right\Vert _{2}$,
with it being no longer the case as in {[}5{]} that $A$ and $P$
commute when acting on general vectors in $L^{2}(\Omega)$. But we
observe that $\left\Vert A^{-\beta}P\left(u_{n}\cdot\nabla\right)u_{n}\right\Vert _{2}=\left\Vert (A^{-\beta}PA^{\beta})A^{-\beta}\left(u_{n}\cdot\nabla\right)u_{n}\right\Vert _{2}$
and that, by using duality arguments similar to those used in {[}17{]},
we have that $A^{-\beta}PA^{\beta}$ is a bounded operator on $L^{2}(\Omega)$,
hence there exists a constant $K_{\beta}$ such that $\left\Vert A^{-\beta}P\left(u_{n}\cdot\nabla\right)u_{n}\right\Vert _{2}\leq K_{\beta}\left\Vert A^{-\beta}\left(u_{n}\cdot\nabla\right)u_{n}\right\Vert _{2}$
and the analysis can now proceed as in {[}5{]} just by incorporating
the additional constant $K_{\beta}$. While this needs to be done
in a number of places, we see that the modifications are quickly performed
in each case, resulting in the proof of Theorem 1. 

For the proof of Theorem 3 we will need higher-order a priori estimates
depending only on the assumed bound given by (1.7). Such bounds are
easily obtained, and the basic method is illustrated in its simplest
form in the following: taking the inner product of both sides of (1.1a)
with $u$, noting that $\left(P\left(u\cdot\nabla\right)u,u\right)=0$
and observing that $(\mu A_{\varphi}u,u)\geq(\mu Q_{m}A^{2}u,u)=(\mu Q_{m}^{2}A^{\alpha}u,u)=\mu\Vert Q_{m}A^{\alpha/2}u\Vert_{2}^{2}$
we obtain the inequality 
\begin{equation}
\frac{d}{dt}\left\Vert u\right\Vert _{2}^{2}+\mu\Vert Q_{m}A^{\alpha/2}u\Vert_{2}^{2}+\nu\left\Vert A^{1/2}u\right\Vert _{2}^{2}\leq\frac{1}{\nu\lambda_{1}}\left\Vert f\right\Vert _{2}^{2}\tag{2.15}
\end{equation}
where we have used in standard fashion Young's inequality and Poincaré's
inequality for the term $\left(f,u\right)$. The term $\mu\Vert Q_{m}A^{\alpha/2}u\Vert_{2}^{2}$
can be discarded from the left-hand side of (2.15) to obtain 
\begin{equation}
\frac{d}{dt}\left\Vert u\right\Vert _{2}^{2}+\nu\left\Vert A^{1/2}u\right\Vert _{2}^{2}\leq\frac{1}{\nu\lambda_{1}}\left\Vert f\right\Vert _{2}^{2}\tag{2.16}
\end{equation}
which is the same energy inequality satisfied by solutions of the
NSE; from this (1.8) follows for $v=u$. 

Using arguments similar to those that led to (2.15) higher-order bounds
can be bootstrapped from the bound (1.9); taking the inner-product
of both sides of (1.6) with $A^{\theta}u$, the left-hand side has
the term $\mu\left(A_{\varphi}u,A^{\theta}u\right)$ and using the
inequality $\left(A_{\varphi}u,A^{\theta}u\right)\geq\left(Q_{m}A^{\alpha+\theta}u,u\right)$
noted above in the remarks following (1.5) we obtain the term $\mu\Vert Q_{m}A^{(\alpha+\theta)/2}u\Vert_{2}^{2}$
on the left-hand side which can again be discarded, from which as
in e.g. {[}4{]} an inequality of the form 
\begin{equation}
\frac{d}{dt}\Vert A^{\theta/2}u\Vert_{2}^{2}+\nu\left\Vert A^{(\theta+1)/2}u\right\Vert _{2}^{2}\leq\frac{2}{\nu}M_{\theta}\Vert A^{\theta/2}u\Vert_{2}^{2}\Vert A^{\theta/2}u\Vert_{2}^{2}+\frac{2}{\nu}\Vert A^{(\theta-1)/2}f\Vert_{2}^{2}\tag{2.17}
\end{equation}
can be derived satisfied by the solutions of both (1.6) and the NSE,
where $M_{\theta}$ depends on $\theta$, the constants $M_{0}$ and
$M_{1}$ appearing in (2.2), and other generic contstants appearing
in typical calculations for the NSE ((see e.g. {[}6{]}). From bounds
on $\Vert A^{\theta/2}u\Vert_{2}$ bounds on $\Vert A^{(\theta+1)/2}u\Vert_{2}$
can be bootstrapped from the inequality (2.17) as demonstrated in
e.g. {[}6{]} to obtain the following result involving uniform higher-order
bounds on each interval $[\tau,T]$ with $0<\tau$ if $f$ in (1.6)
and the NSE is smooth enough. 

\begin{theorem} Let $[0,T]$ be any interval on which $f\in C([0,T];D(A^{(\theta-1)/2})$)
for a natural number $\theta\geq2$. Then given the uniform $H^{1}$-bound
(1.9) for solutions $v=u$ of (1.6) or and the NSE on $[0,T]$, there
exists on each subinterval $[\tau,T]$ with $\tau>0$ a constant $U_{\theta,\tau}$
such that $\sup_{\tau\leq t\leq T}\Vert A^{\theta/2}u(t)\Vert_{2}\leq U_{\theta,\tau}$
and a constant $U_{\theta,\tau}^{'}$ such that $\int_{0}^{T}\Vert A^{(\theta+1)/2}u(t)\Vert_{2}\leq U_{\theta,\tau}^{'}$The
constants $U_{\theta,\tau}$ and $U_{\theta,\tau}^{'}$ depend only
on $\tau$ and $U_{1}$. \end{theorem}

In Theorem 4 it is assumed that $\tau>0$ since it is standard to
assume that $u_{0}$ has no more than $H^{1}$-regularity. But by
using the regularity of solutions for $t>0$ we can by replacing $u(x,0)$
by $u(x,\tau)$ if necessary assume in what follows that $\tau=0$
for simplicity, in which case we can replace $U_{\theta,\tau}$ by
$U_{\theta}$ and $U_{\theta,\tau}^{'}$ by $U_{\theta}^{'}$ in Theorem
3 for each natural number $\theta\geq2$. 

\section{Proof of Theorems 2 and 3}

We first prove Theorem 3; let $u$ be the solution of (1.6), let $v$
be the solution of (1.2), and let $w=u-v$. We assume for simplicity
that $u$ and $v$ share the same data; modifying accordingly to obtain
a result allowing a continuous dependence on data will be seen to
be sraightforward. To prove Theorem 3, we subtract the NSE from (1.6)
and take the inner-product of both sides of the resulting equation
for $w$ with $A^{\theta}w$; using the inequality $\left(A_{\varphi}u,A^{\theta}u\right)\geq\left(Q_{m}A^{\alpha+\theta}u,u\right)$
noted above we see that $\left(A_{\varphi}u,A^{\theta}u\right)\geq(\mu Q_{m}A^{\alpha+\theta-1}u,Aw)=\left\Vert Q_{m}A^{(\alpha+\theta-1)/2}u\right\Vert _{2}^{2}$
and treating the term on the right-hand side as a forcing term when
using Young's inequality we use calculations similar to those leading
to (2.17) to obtain that 
\begin{equation}
\frac{d}{dt}\Vert A^{\theta/2}w\Vert_{2}^{2}+\nu\left\Vert A^{(\theta+1)/2}w\right\Vert _{2}^{2}\leq\frac{2}{\nu}M_{\theta}\Vert A^{\theta/2}w\Vert_{2}^{2}[\Vert A^{\theta/2}u\Vert_{2}^{2}+\Vert A^{\theta/2}v\Vert_{2}^{2}]+\frac{2\mu}{\nu}\left\Vert Q_{m}A^{(\alpha+\theta-1)/2}u\right\Vert _{2}^{2}.\tag{3.1}
\end{equation}
Integrating both sides of (3.1) we have that 
\begin{equation}
\Vert A^{\theta/2}w\Vert_{2}^{2}\leq\mu\int_{0}^{T}\frac{2}{\nu}\left\Vert Q_{m}A^{(\alpha+\theta-1)/2}u\right\Vert _{2}^{2}ds+\int_{0}^{t}\frac{2}{\nu}M_{\theta}[\Vert A^{\theta/2}u\Vert_{2}^{2}+\Vert A^{\theta/2}v\Vert_{2}^{2}]\Vert A^{\theta/2}w\Vert_{2}^{2}ds\tag{3.2}
\end{equation}
where we have used that $w(0)=0$ and that $\int_{0}^{t}\frac{2}{\nu}\left\Vert Q_{m}A^{(\alpha+\theta-1)/2}u\right\Vert _{2}^{2}ds\leq\int_{0}^{T}\frac{2}{\nu}\left\Vert Q_{m}A^{(\alpha+\theta-1)/2}u\right\Vert _{2}^{2}ds$.
Applying Theorem 3 and the remarks following to (3.2) as well as Gronwall's
inequality we obtain that 
\begin{equation}
\Vert A^{\theta/2}w\Vert_{2}^{2}\leq\mu\left[\int_{0}^{T}\frac{2}{\nu}\left\Vert Q_{m}A^{(\theta+3)/2}u\right\Vert _{2}^{2}ds\right]\exp\left(\frac{2}{\nu}M_{\theta}(U_{\theta}^{'}+V_{\theta}^{'})\right)\tag{3.3}
\end{equation}
where $V_{\theta}$ and $V_{\theta}^{'}$ are the same as $U_{\theta}$
and $U_{\theta}^{'}$ but for the solution $v$ of (1.2). For the
case $\mu\rightarrow0$, the assumptions on $f$ imply that Theorem
3 holds with $\theta$ replaced by $\alpha+\theta-1$, and then we
simply note that $\int_{0}^{T}\frac{2}{\nu}\left\Vert Q_{m}A^{(\alpha+\theta-1)/2}u\right\Vert _{2}^{2}ds\leq\frac{2}{\nu}U_{\alpha+\theta-2}^{'}$.
For the case $m\rightarrow\infty$ we note that since $\left\Vert Q_{m}A^{(\alpha+\theta-1)/2}u\right\Vert _{2}^{2}\leq\left\Vert A^{(\alpha+\theta-1)/2}u\right\Vert _{2}^{2}$
we have from Theorem 3 that $\left\Vert Q_{m}A^{(\alpha+\theta-1)/2}u(t)\right\Vert _{2}^{2}\rightarrow0$
as $m\rightarrow\infty$ for a.e. $t$ and hence that $\int_{0}^{T}\left\Vert Q_{m}A^{(\alpha+\theta-1)/2}u\right\Vert _{2}^{2}ds\rightarrow0$
as $m\rightarrow\infty$ by the dominated convergence theorem, thus
establishing Theorem 3.

We now prove Theorem 2; this requires a more complex semigroup approach,
since we want to only rely on (1.9). By the ADC assumptions on $A_{\varphi}$
and the functional calculus we have that $\mu A_{\varphi}$ generates
a contraction semigroup as does $\nu A$, and since these two operators
commute we have that $\exp(\nu A+\mu A_{\varphi})=\exp\nu A\exp\mu A_{\varphi}$;
for solutions $u$ of (1.1) we thus have that $u$ satisfies 
\begin{equation}
u(t)=e^{-t(\nu A+\mu A_{\varphi})}u_{0}+\int_{0}^{t}e^{-(t-s)(\nu A+\mu A_{\varphi})}P\left(u\cdot\nabla\right)uds+\int_{0}^{t}e^{-(t-s)(\nu A+\mu A_{\varphi})}fds\tag{3.4}
\end{equation}
for $t\in[0,T]$, while for solutions $v$ of the NSE we have that
$v$ satisfies (3.4) with $\mu=0$; subtracting, we have for $w=u-v$
that 
\begin{equation}
w(t)=h(t)+e^{-t(\nu A+\mu A_{\varphi})}u_{0}+\int_{0}^{t}e^{-(t-s)\nu A}e^{-(t-s)\mu A_{\varphi}}P(\left(w\cdot\nabla\right)u+\left(v\cdot\nabla\right)w)ds\tag{3.5}
\end{equation}
wherein
\begin{align}
h(t)\equiv e^{-t\nu A}(e^{-t\mu A_{\varphi}}-I)u_{0}+\int_{0}^{t}e^{-(t-s)\nu A}(e^{-(t-s)\mu A_{\varphi}}-I)fds\nonumber \\
+\int_{0}^{t}e^{-(t-s)\nu A}(e^{-(t-s)\mu A_{\varphi}}-I)P(\left(w\cdot\nabla\right)u+\left(v\cdot\nabla\right)w)ds.\tag{3.6}
\end{align}
Before analyzing (3.5), (3.6) we note that using (2.2) we have for
a constant $M_{1}$ that $\Vert A^{-1/4}P(v\cdot\nabla)w\Vert_{2}\leq M_{0}\Vert P(u\cdot\nabla)u\Vert_{3/2}$,
and we recall that $P$ is a bounded operator on $L^{p}(\Omega)$,
$1<p<\infty$ (see e.g. {[}16{]}, {[}17{]} and the references contained
therein); absorbing its operator bound for $p=3/2$ into $M_{1}$
and using Hölder's inequality we thus have that $\Vert A^{-1/4}P(v\cdot\nabla)w\Vert_{2}\leq M_{1}\Vert v\Vert_{6}\Vert\nabla w\Vert_{2}$,
from which it again follows from (2.2) and a suitably redefined $M_{1}$
that ies 
\begin{equation}
\Vert A^{-1/4}P(v\cdot\nabla)w\Vert_{2}\leq M_{1}\Vert A^{1/2}v\Vert_{2}\Vert A^{1/2}w\Vert_{2}.\tag{3.7}
\end{equation}
We also recall that $e^{-t\nu A}$ is an analytic semigroup (see e.g.
{[}16{]}, {[}17{]} and the references contained therein), thus there
is a constant $c$ such that for each $\gamma,b>0$ 
\begin{equation}
\Vert A^{\gamma}e^{-tbA}u\Vert_{2}\leq\frac{c}{(tb)^{\gamma}}\Vert u\Vert_{2}.\tag{3.8}
\end{equation}
We now replace $t$ with $\tau$ in (3.5) with $0\leq\tau\leq t$
and apply $A^{1/2}$ to both sides to obtain
\begin{equation}
\Vert A^{1/2}w(\tau)\Vert_{2}\leq\Vert A^{1/2}h(\tau)\Vert_{2}+\int_{0}^{\tau}\Vert A^{3/4}e^{-(\tau-s)\nu A}A^{-1/4}P(\left(w\cdot\nabla\right)u+\left(v\cdot\nabla\right)w)\Vert_{2}ds\tag{3.9}
\end{equation}
and in (3.9) we have used the fact that $e^{-(\tau-s)\nu A}$ and
$e^{-(\tau-s)\mu A_{\varphi}}$ are contraction semigroups; we will
use (3.7) to show that the terms on the right-hand side of (3.9) are
well-defined in terms of our assumed estimates. For a function $g\in H$
and $n\geq m$, using the factorization $e^{-(t-s)\nu A}=e^{-[(t-s)/2]\nu A}e^{-[(t-s)/2]\nu A}$
as well as (3.7) and (3.8) we have that
\begin{align}
 & \Vert A^{3/4}e^{-(t-s)\nu A}g\Vert_{2}\leq\Vert P_{n}A^{3/4}e^{-(t-s)\nu A}g\Vert_{2}+\Vert Q_{n}A^{3/4}e^{-(t-s)\nu A}g\Vert_{2}\nonumber \\
 & \leq\Vert A^{3/4}P_{n}g\Vert_{2}+\Vert A^{3/4}e^{-[(t-s)/2]\nu A}\Vert_{2}\Vert e^{-[(t-s)/2]\nu A}Q_{n}g\Vert_{2}\nonumber \\
 & \leq\lambda_{n}^{3/4}\Vert P_{n}g\Vert_{2}+\frac{2c}{\nu(t-s)^{3/4}}e^{-(\nu/2)\lambda_{n+1}(t-s)}\Vert Q_{n}g\Vert_{2}\tag{3.10}
\end{align}
since $\lambda_{n}$ is the largest eigenvalue of $A$ on $P_{n}H$
and $\lambda_{n+1}$ is the smallest eigenvalue of $A$ on $Q_{n}H$.
Applying (3.10) to $g=A^{-1/4}P(\left(w\cdot\nabla\right)u+\left(v\cdot\nabla\right)w)$
in (3.10), noting that $P_{n}$ and $Q_{n}$ are projections, and
using (3.6) we obtain that
\begin{align}
\Vert A^{1/2}w(\tau)\Vert_{2}\leq\Vert A^{1/2}h(\tau)\Vert_{2}+\int_{0}^{\tau}\lambda_{n}^{3/4}M_{1}[\Vert A^{1/2}u\Vert_{2}+\Vert A^{1/2}v\Vert_{2}]\Vert A^{1/2}w\Vert_{2}ds\nonumber \\
+\int_{0}^{\tau}\frac{2c}{\nu(\tau-s)^{3/4}}e^{-(\nu/2)\lambda_{n+1}(\tau-s)}M_{1}[\Vert A^{1/2}u\Vert_{2}+\Vert A^{1/2}v\Vert_{2}]\Vert A^{1/2}w\Vert_{2}ds.\tag{3.11}
\end{align}
Applying (1.9) to (3.11), setting $\rho(t)\equiv\sup_{0\leq s\leq t}\Vert A^{1/2}w(s)\Vert_{2}$,
and changing variables in the second integral we have for $G(\tau)\equiv\Vert A^{1/2}h(\tau)\Vert_{2}$
that
\begin{align}
\Vert A^{1/2}w(\tau)\Vert_{2}\leq G(\tau)+\int_{0}^{\tau}\lambda_{n}^{3/4}M_{1}[\Vert A^{1/2}u\Vert_{2}+\Vert A^{1/2}v\Vert_{2}]\Vert A^{1/2}w\Vert_{2}ds\nonumber \\
+\frac{4c}{\nu}M_{1}U_{1}\rho(t)\int_{0}^{\tau}\frac{e^{-(\nu/2)\lambda_{n+1}s}}{s^{3/4}}ds.\tag{3.12}
\end{align}
Replacing $\tau$ by $t$ on the right-hand-side of (3.12) and using
the estimate $\int_{0}^{t}\frac{e^{-\beta s}}{s^{\gamma}}ds\leq\frac{1}{1-\gamma}\frac{e^{-\gamma}}{\beta^{1-\gamma}}\leq\frac{1}{(1-\gamma)\beta^{1-\gamma}}$
(see e.g. {[}1{]} or {[}14{]}) for $\beta=(\nu/2)\lambda_{n+1}$ and
$\gamma=3/4$ we have for $\overline{G}(T)\equiv\sup_{0\leq\tau\leq T}G(\tau)$
that
\begin{equation}
\Vert A^{1/2}w(\tau)\Vert_{2}\leq\overline{G}(T)+\int_{0}^{t}\lambda_{n}^{3/4}M_{1}[\Vert A^{1/2}u\Vert_{2}+\Vert A^{1/2}v\Vert_{2}]\Vert A^{1/2}w\Vert_{2}ds+\frac{16c}{\nu}M_{1}U_{1}\frac{1}{[(\nu/2)\lambda_{n+1}]^{1/4}}\rho(t)\tag{3.13}
\end{equation}
and since (3.13) holds for all $0\leq\tau\leq t$ on the left-hand
side we have that
\begin{equation}
\rho(t)\leq\overline{G}(T)+\int_{0}^{t}\lambda_{n}^{3/4}M_{1}[\Vert A^{1/2}u\Vert_{2}+\Vert A^{1/2}v\Vert_{2}]\Vert A^{1/2}w\Vert_{2}ds+\frac{16c}{\nu}M_{1}U_{1}\frac{1}{[(\nu/2)\lambda_{n+1}]^{1/4}}\rho(t).\tag{3.14}
\end{equation}
We now choose $n\geq m$ large enough so that $16c\nu^{-1}M_{1}U_{1}[(\nu/2)\lambda_{n+1}]^{-1/4}\leq1/2$;
applying this to (3.14), combining terms, and replacing $\Vert A^{1/2}w(s)\Vert_{2}$
by its upper bound $\rho(s)=\sup_{0\leq\eta\leq s}\Vert A^{1/2}w(\eta)\Vert_{2}$
we have that
\begin{equation}
\rho(t)\leq2\overline{G}(T)+\int_{0}^{t}4\lambda_{n}^{3/4}M_{1}[\Vert A^{1/2}u\Vert_{2}+\Vert A^{1/2}v\Vert_{2}]\rho(s)ds.\tag{3.15}
\end{equation}
To apply Gronwall's inequality to (3.15) we first need to estimate
$\overline{G}(T)$. In this regard we first establish that $G(\tau)$
as defined in the line following (3.11) is continuous in $\tau$.
This is fairly clear but for completeness we sketch the details. Defining
$G(t)$ by (3.9) with $\tau$ replaced by $t$, using the fact that
$|\Vert\phi\Vert_{2}-\Vert\psi\Vert_{2}|\leq\Vert\phi-\psi\Vert_{2}$
for suitable vectors $\phi$ and $\psi$, and using properties of
semigroups as well as (3.8) we have for $\phi=f+P\left(v\cdot\nabla\right)v$
that 
\begin{align}
 & |G(t)-G(\tau)|\leq\Vert e^{-(t-\tau)\nu A}(e^{-(t-\tau)\mu A_{\varphi}}-I)u_{0}\Vert_{2}+\int_{\tau}^{t}\Vert A^{3/4}e^{-(t-s)\nu A}(e^{-(t-s)\mu A_{\varphi}}-I)A^{-1/4}\phi\Vert_{2}ds\nonumber \\
 & +\int_{0}^{\tau}\Vert A^{3/4}e^{-(\tau-s)\nu A}[(e^{-(t-\tau)\nu A}-I)+e^{-(\tau-s)\mu A_{\varphi}}(e^{-(t-\tau)\nu A}e^{-(t-\tau)\mu A_{\varphi}}-I)]A^{-1/4}\phi\Vert_{2}ds\nonumber \\
 & \leq\Vert(e^{-(t-\tau)\mu A_{\varphi}}-I)u_{0}\Vert_{2}+\int_{\tau}^{t}\frac{c}{\nu(t-s)^{3/4}}\Vert(e^{-(t-s)\mu A_{\varphi}}-I)A^{-1/4}\phi\Vert_{2}ds\nonumber \\
 & +\int_{0}^{\tau}\frac{c}{\nu(\tau-s)^{3/4}}[\Vert(e^{-(t-\tau)\nu A}-I)A^{-1/4}\phi\Vert_{2}+\Vert(e^{-(t-\tau)(\nu A+\mu A_{\varphi})}-I)A^{-1/4}\phi\Vert_{2}]ds.\tag{3.16}
\end{align}
Now $A^{-1/4}\phi\in L^{2}(\Omega)$ by (1.7), (3.7) and the assumptions
on $f$, thus $\Vert(e^{-(t-\tau)\mu A_{\varphi}}-I)u_{0}\Vert_{2}\rightarrow0$
as $\tau\rightarrow t$ and for each $s$ $\Vert(e^{(t-\tau)B}-I)A^{-1/4}\phi(s)\Vert_{2}\rightarrow0$
as $\tau\rightarrow t$ for $B=-\nu A,\mu A_{\varphi},-(\nu A+\mu A_{\varphi})$
since each of these operators generates a strongly-continuous semigroup.
The integrand in the third line of (3.16) is bounded by $2c\nu^{-1}(t-s)^{-3/4}\Vert A^{-1/4}\phi\Vert_{2}\in L^{1}(0,t))$
and so the corresponding integral goes to zero as $\tau\rightarrow t$.
Similar arguments for the integrands in the fourth line of (3.16)
together with the above observations on strong convergence shows that
the corresponding integrals go to zero as $\tau\rightarrow t$ by
the dominated convergence theorem. Thus from (3.16) we have that $G(\tau)$
is continuous on $[0,T]$, hence there is a $t_{0}\in[0,T]$ such
that $\overline{G}(T)=G(t_{0}).$ Meanwhile from the inequality $ab\leq\frac{1}{2}(a^{2}+b^{2})$
with $b=1$ we have that $\Vert A^{1/2}u\Vert_{2}+\Vert A^{1/2}v\Vert_{2}\leq\frac{1}{2}(\Vert A^{1/2}u\Vert_{2}^{2}+\Vert A^{1/2}v\Vert_{2}^{2})+1$
and hence from (2.4), (3.15) and Gronwall's inequality we have that
\begin{equation}
\rho(t)\leq2G(t_{0})\exp[2\lambda_{n}^{3/4}M_{1}(\frac{1}{\nu^{2}\lambda_{1}}\int_{0}^{T}\left\Vert f\right\Vert _{2}^{2}ds+2T)].\tag{3.17}
\end{equation}
Using (3.8) and arguments similar to those used in (3.10) we have
for each $t\in[0,T]$ that
\begin{align}
\Vert A^{1/2}h(t)\Vert_{2}\leq\parallel(e^{-t\mu A_{\varphi}}-I)u_{0}\parallel_{2}+\int_{0}^{t}\parallel(e^{-(t-s)\mu A_{\varphi}}-I)A^{1/2}f\parallel_{2}ds\nonumber \\
+\int_{0}^{t}\frac{c}{\nu(t-s)^{3/4}}\parallel(e^{-(t-s)\mu A_{\varphi}}-I)A^{-1/4}P(\left(w\cdot\nabla\right)u+\left(v\cdot\nabla\right)w))\parallel_{2}ds.\tag{3.18}
\end{align}
For each $\tau\in[0,t]$ if we now let $\mu\rightarrow0$ then $\parallel(e^{-\tau\mu A_{\varphi}}-I)\phi\parallel_{2}\rightarrow0$
for each $\phi\in L^{2}(\Omega)$ so by arguments similar to those
above in analyzing (3.16) we have that $\Vert A^{1/2}h(t)\Vert_{2}\rightarrow0$
as $\mu\rightarrow0$ for each $t$. Setting $t=t_{0}$ we have that
$G(t_{0})\rightarrow0$ as $\mu\rightarrow0$ and hence that $\Vert A^{1/2}(u(t)-v(t)\Vert_{2}\rightarrow0$
uniformly for $t\in[0,T]$ as $\mu\rightarrow0$ follows from (3.18).
For the case $m_{0}\rightarrow\infty$ we note that $A_{\varphi}=Q_{m_{0}}A_{\varphi}$
so that by the functional calculus $\parallel(e^{-\tau\mu A_{\varphi}}-I)\phi\parallel_{2}=\parallel Q_{m_{0}}(e^{-\tau\mu A_{\varphi}}-I)\phi\parallel_{2}$
for each $\tau\in[0,t]$ and each $\phi\in L^{2}(\Omega)$. Since
$\parallel Q_{m_{0}}\psi\parallel_{2}\rightarrow0$ as $m_{0}\rightarrow\infty$
for each $\psi\in L^{2}(\Omega)$, setting $\psi=(e^{-\tau\mu A_{\varphi}}-I)\phi$
for each $\tau$ we have again by similar arguments to those used
above that $G(t_{0})\rightarrow0$ but now as $m_{0}\rightarrow\infty$.
It thus follows that $\Vert A^{1/2}(u(t)-v(t)\Vert_{2}\rightarrow0$
uniformly for $t\in[0,T]$ as $m_{0}\rightarrow\infty$, and with
this we complete the proof of Theorem 2. 

\section{Reformulation of the HNSE and SHNSE}

We begin by considering (1.5); using (1.4) and the formula $Au=-\Delta u+\nabla p_{s}(u)$
noted above in the remarks preceeding (1.2) we have that 
\begin{equation}
A^{2}u=P(-\Delta)^{2}u=P(-\Delta)(-\Delta)u=(-\Delta)P(-\Delta)u+\nabla p_{s}((-\Delta)u)\tag{4.1}
\end{equation}
but in turn $P(-\Delta)u=(-\Delta)Pu+\nabla p_{s}(u)=(-\Delta)u+\nabla p_{s}(u)$,
so combining with (4.1) we have that 
\begin{align}
A^{2}u & =(-\Delta)[(-\Delta)u+\nabla p_{s}(u)]+\nabla p_{s}((-\Delta)u)\nonumber \\
 & =(-\Delta)^{2}u+(-\Delta)(\nabla p_{s}(u))+\nabla p_{s}((-\Delta)u)\nonumber \\
 & =(-\Delta)^{2}u+\nabla p_{s}((-\Delta)u)\tag{4.2}
\end{align}
 where similarly to the remarks following (1.3) we see that $(-\Delta)(\nabla p_{s}(u))=\nabla(-\Delta p_{s}(u))=0$,
and where we note that $\nabla p_{s}((-\Delta)u)$ satisfies
\begin{align}
\Delta p_{s}((-\Delta)u)=0,x\in\Omega,\tag{4.3a}\nonumber \\
\mathbf{n}\cdot\nabla p_{s}((-\Delta)u)=-\mathbf{n\cdot}(-\Delta)^{2}u,x\in\Gamma\tag{4.3b}
\end{align}
which is obtained by substituting $(-\Delta)u$ in place of $u$ in
(1.2). With the formula $PA^{2}u=(-\Delta)^{2}u+\nabla p_{s}((-\Delta)u)$
where $\nabla p_{s}((-\Delta)u)$ satisfies (4.3), we have, by following
the Stokes-pressure methodology of {[}29{]}, {[}30{]}, retained its
philosophy of replacing the dependence on $P$ and the Stokes operator
with a relatively-simple operator and a relatively-simple elliptic
boundary-value problem for potential computational applications. Using
the linearity of the systems (1.2) and (4.3) and again using that
$Au=P(-\Delta)u=(-\Delta)u+\nabla p_{s}(u)$ we have that (1.5) can
be reformulated as the system
\begin{align}
\frac{d}{dt}u+\mu(-\Delta)^{2}u+\nu(-\Delta)u+P\left(u\cdot\nabla\right)u+\nabla p_{s}(\mu(-\Delta)u+\nu u)=f,\tag{4.4a}\nonumber \\
u(x,0)=u_{0}(x)\tag{4.4b}
\end{align}
where $\mu\nabla p_{s}((-\Delta)u)+\nu\nabla p_{s}(u)=\nabla p_{s}(\mu(-\Delta)u+\nu u)$
satisfies
\begin{align}
\Delta p_{s}(\mu(-\Delta)u+\nu u)=0,x\in\Omega,\tag{4.5a}\nonumber \\
\mathbf{n}\cdot\nabla p_{s}(\mu(-\Delta)u+\nu u)=\mathbf{n\cdot}\Delta(\nu I-\mu\Delta)u,x\in\Gamma.\tag{4.5b}
\end{align}
Meanwhile the operator $(-\Delta)^{2}u$ is relatively easy to computationally
implement in the interior of the domain and the boundary condition
(4.5b) is relatively straighforward to implement as well. To implement
(1.6), we have that $\mu(-\Delta)^{2}u+\nabla p_{s}(\mu(-\Delta)u)=\mu[(-\Delta)^{2}u+\nabla p_{s}((-\Delta)u)]=\mu P_{n}^{*}[(-\Delta)^{2}u+\nabla p_{s}((-\Delta)u)]+\mu Q_{n}^{*}[(-\Delta)^{2}u+\nabla p_{s}((-\Delta)u)]$
where $P_{n}^{*}$ is the projection onto the first $n$ eigenspaces
of $A$ or onto suitable psuedospectral basis spaces as may arise
in computational settings, and $Q_{n}^{*}=I-P_{n}^{*}$. Again arguing
that the effects of $\mu P_{n}^{*}[(-\Delta)^{2}u+\nabla p_{s}((-\Delta)u)]$
are small on the elements corresponding to lower frequencies below
some cutoff $m$, we replace $\mu(-\Delta)^{2}u+\mu\nabla p_{s}((-\Delta)u)$
in (4.4a) by $\mu Q_{n}^{*}[(-\Delta)^{2}u+\nabla p_{s}((-\Delta)u)]$
for $n\geq m$ which in light of (4.2) is a version of a special case
of (1.6) that represents an alternative formulation for computations. 

Now for each smooth vector $h\in L^{2}(\Omega)$ we have the Helmholtz
decomposition $h=Ph+\nabla\phi(h)$ where $\Delta\phi=\nabla\cdot h$
in $\Omega$ and $\mathbf{n}\cdot\nabla\phi(h)=-\mathbf{n\cdot}h$
on $\Gamma$. Setting $h=g-(u\cdot\nabla)u$, assuming as in {[}29{]},
{[}30{]} that the pressure is normalized to mean zero,and comparing
(4.4a) with (1.6a) for $A_{\varphi}=P(-\Delta)^{2}$ we have that
\begin{equation}
p=p_{s}(\nu u+\mu(-\Delta)u)+\phi(g-(u\cdot\nabla)u)\tag{4.6}
\end{equation}
or alternatively, for $Q=I-P$, 
\begin{equation}
\nabla p=\nabla p_{s}(\nu u+\mu(-\Delta)u)+Q(g-(u\cdot\nabla)u)\tag{4.7}
\end{equation}
so that as in {[}32{]}, {[}33{]} we have an expression for $p$ that
is explicit in terms of solutions of boundary value problems given
(4.5) and the definition of $\phi$.

We can generalize these developments for the purposes of robustness
as in {[}29{]}, {[}30{]} to embed (4.4), (4.5) into a larger class
of equations in which the divergence-free condition no longer necessarilty
holds. In this case from the development following (4.1) we have $P(-\Delta)u=(-\Delta)Pu+\nabla p_{s}(u)$
where we no longer assume that $Pu=u$, in which case as shown in
{[}29{]}, {[}30{]} we have that 
\begin{align}
\Delta p_{s}(u)=0,x\in\Omega,\tag{4.8a}\nonumber \\
\mathbf{n}\cdot\nabla p_{s}(u)=\mathbf{n\cdot(\Delta-\nabla\nabla\cdot)}u,x\in\Gamma\tag{4.8b}
\end{align}
and in which case (4.2) becomes 
\begin{align}
A^{2}u & =(-\Delta)^{2}Pu+\nabla p_{s}((-\Delta)u)\nonumber \\
 & =(-\Delta)^{2}u+(-\Delta)(\Delta(I-P)u)+\nabla p_{s}((-\Delta)u)\nonumber \\
 & =(-\Delta)^{2}u+(-\Delta)\nabla\nabla\cdot u+\nabla p_{s}((-\Delta)u)\tag{4.9}
\end{align}
where we have noted that $\Delta(I-P)u=\Delta\nabla\phi(u)=\nabla\Delta\phi=\nabla\nabla\cdot u$
using the definition of $\phi$. Since also $(-\Delta)\nabla\nabla\cdot u=\nabla\nabla\cdot(-\Delta u)$,
we have from (4.9) the following system as a generalized version of
(4.4)
\begin{align}
\frac{d}{dt}u+\mu(-\Delta)^{2}u+\mu\nabla\nabla\cdot(-\Delta)u+\nu(-\Delta)u+P\left(u\cdot\nabla\right)u+\nabla p_{s}(\mu(-\Delta)u+\nu u)=f,\tag{4.10a}\nonumber \\
u(x,0)=u_{0}(x),\tag{4.10b}
\end{align}
where $\nabla p_{s}(\mu(-\Delta)u+\nu u)$ now satisfies
\begin{align}
\Delta p_{s}(\mu(-\Delta)u+\nu u)=0,x\in\Omega,\tag{4.11a}\nonumber \\
\mathbf{n}\cdot\nabla p_{s}(\mu(-\Delta)u+\nu u)=\mathbf{n\cdot}(\Delta-\nabla\nabla\cdot)(\nu I-\mu\Delta)u,x\in\Gamma,\tag{4.11b}
\end{align}
and the pressure $p$ now satisfies 
\begin{equation}
p=p_{s}(\nu u+\mu(-\Delta)u)+\mu\nabla\cdot(-\Delta)u+\phi(g-(u\cdot\nabla)u).\tag{4.12}
\end{equation}
Now since $\nabla\cdot(\mu\nabla\nabla\cdot(-\Delta)u)=\mu\Delta\nabla\cdot(-\Delta)u=-\mu\nabla\cdot(-\Delta)^{2}u$
we have by taking the divergence of both sides of (4.10a) and using
(4.10), (4.11), and the definition of $\phi$ that $\nabla\cdot u$
satisthat 
\begin{align}
(\nabla\cdot u)_{t}=\nu\Delta(\nabla\cdot u),x\in\Omega,\tag{4.13a}\nonumber \\
\mathbf{n}\cdot\nabla(\nabla\cdot u)=0,x\in\varGamma.\tag{4.13b}
\end{align}
Thus as for the equation in which the NSE was embedded in {[}29{]},
{[}30{]} we have in particular that $(\nabla\cdot u)=0$ initially
if and only if $(\nabla\cdot u)=0$ for all later time and more generally
that if $\nabla\cdot u$ is small initially then $(\nabla\cdot u)=0$
remains small for all time by the maximum principle. Hence the system
(4.10), (4.11) extends the (constrained) dynamics of (4.4), (4.5)
in a well-posed manner as was also the case for the extended version
of the NSE in {[}32{]}, {[}33{]}. We can generalize (1.6) to this
setting by arguing similarly as before and replacing $\mu[(-\Delta)^{2}u+\mu\nabla\nabla\cdot(-\Delta)u+\nabla p_{s}((-\Delta)u)]$
in (4.10a) by $\mu Q_{n}^{*}[(-\Delta)^{2}u+\mu\nabla\nabla\cdot(-\Delta)u+\nabla p_{s}((-\Delta)u)]$
for suitable $n$; since this is equivalent to (4.10) with $f$ replaced
by $f_{n}^{*}=f-g_{n}^{*}$ where $g_{n}^{*}=\mu P_{n}^{*}[(-\Delta)^{2}u+\mu\nabla\nabla\cdot(-\Delta)u+\nabla p_{s}((-\Delta)u)]$;
robustness is preserved by continuous dependence on data in (4.10),
and (4.13) is modified by adding $-g_{n}^{*}$ to the right-hand side
of (4.13a); robustness is preserved by continuous dependence on data
on intervals $[0,T]$ given that $g_{n}^{*}$ is assumed to be negligibly
close to zero. With these considerations our discussion of alternative
formulations for computational implementation of (1.5) and (1.6) is
complete.

\section{Conclusion}

The proof of Theorem 2 necessarily is broken into high- and low-frequency
components by the need to absorb powers of $A$ using (3.6). With
decomposition as in (3.10) and (3.17) Gronwall's inequality can be
used for the $P_{n}$-terms since they can absorb powers of $A$ as
bounded operators, but (3.6) prevents the use of Gronwall's inequality
for the $Q_{n}$-terms. Instead, we are able to use (3.6) to handle
the $Q_{n}$-terms due to our ability to make these terms small and
absorb them for large enough $n$. The use of (3.6) and the underlying
use of the estimate (3.8) distinguishes these arguments from the spectral-decomposition
methods used in {[}4{]}. We also note that the only place in the proof
that the assumption (1.9) is used is in the selection of $n$ for
the absorbing of the $Q_{n}$-terms, otherwise only the standard energy
estimate (2.4) is needed. 

As noteed in the introduction the results and discussion in {[}25{]},
{[}27{]} show that SGS dissipation can potentially under time averaging
be a reasonably good model of the effects of the unresolved inertial-range
scales on the resolved scales. This suggests in particular that given
their independence of time-averaging steady-state solutions of (1.1)
as well as trajectories close to them would be worthy topics for further
study. 

Besides the motivations discussed in {[}5{]}, {[}18{]}, {[}19{]} the
potential of the SHNSE to serve as a significant and reasonable LES
turbulence model is suggested further by considering the successful
implementations of spectral vanishing viscosity (SVV) in modeling
turbulence as noted in the introduction. Indeed SVV can be seen as
an implementation of the SHNSE in some sense as noted in {[}21{]}
and resembles the SHNSE in truncation as noted in {[}5{]}. In {[}31{]}
SVV is used in particular in a high-Reynolds-number wind-tunnel simulation
in which good results are obtained overall and any significant deviation
from expected results seems to be localized at the boundary. Other
potential numerical issues for the SHNSE involving the boundary were
identified in {[}18{]}. These observations suggest that boundary issues
represent a central focus in modeling accuracy in computational implementations
of the SVV and SHNSE that reflect their importance in theoretical
issues as discussed here. While the formulation for the SHNSE on bounded
domains discussed here may not directly impact these computational
issues, it may serve as a foundational and conceptual framework from
which to approach successful implementations in the future.

\end{document}